\documentclass[12pt]{article}

\usepackage[colorlinks,allcolors=blue]{hyperref}

\usepackage[margin=1in]{geometry}
\usepackage{amsfonts, amsmath, amssymb, amsthm} 
\usepackage[none]{hyphenat}
\usepackage{fancyhdr} 
\usepackage[utf8]{inputenc}
\usepackage{mathtools}

\makeatletter
\def\legendre@dash#1#2{\hb@xt@#1{%
  \kern-#2\p@
  \cleaders\hbox{\kern.5\p@
    \vrule\@height.2\p@\@depth.2\p@\@width\p@
    \kern.5\p@}\hfil
  \kern-#2\p@
  }}
\def\@legendre#1#2#3#4#5{\mathopen{}\left(
  \sbox\z@{$\genfrac{}{}{0pt}{#1}{#3#4}{#3#5}$}%
  \dimen@=\wd\z@
  \kern-\p@\vcenter{\box0}\kern-\dimen@\vcenter{\legendre@dash\dimen@{#2}}\kern-\p@
  \right)\mathclose{}}
\newcommand\legendre[2]{\mathchoice
  {\@legendre{0}{1}{}{#1}{#2}}
  {\@legendre{1}{.5}{\vphantom{1}}{#1}{#2}}
  {\@legendre{2}{0}{\vphantom{1}}{#1}{#2}}
  {\@legendre{3}{0}{\vphantom{1}}{#1}{#2}}
}
\def\dlegendre{\@legendre{0}{1}{}}
\def\tlegendre{\@legendre{1}{0.5}{\vphantom{1}}}
\makeatother

\usepackage{booktabs, siunitx,caption}

\newtheorem{theorem}{Theorem}

\title{Primitive Indexes, Zsigmondy Numbers, and Primoverization}
\author{Tejas R. Rao}
\date{September 22, 2018}

\begin{document}
\maketitle

\begin{center}
\small \textbf{Abstract.} We define a \emph{primitive index} of an integer in a sequence to be the index of the term with the integer as a primitive divisor. For the sequences $k^u+h^u$ and $k^u-h^u$, we discern a formula to find the primitive indexes of any composite number given the primitive indexes of its prime factors. We show how this formula reduces to a formula relating the multiplicative order of $k$ modulo $N$ to that of its prime factors. We then introduce immediate consequences of the formula: certain sequences which yield the same primitive indexes for numbers with the same unique prime factors, an expansion of the lifting the exponent lemma for $\nu_2(k^n+h^n)$, a simple formula to find any Zsigmondy number, a note on a certain class of pseudoprimes titled overpseudoprime, and a proof that numbers such as Wagstaff numbers are either overpseudoprime or prime.\\
\end{center}

\setlength{\parindent}{0cm}

\section{The Formula}
A \emph{primitive index} of a number is the index of the term in a sequence that has that number as a primitive divisor. Let $k$ and $h$ be positive integers for the duration of this paper. We denote the primitive index of $N$ in $k^u-h^u$ as $P_N(k,h)$. For $k^u+h^u$, it is denoted $P^N(k,h)$. Note that $P_N(k,1)=O_N(k)$, where $O_N(k)$ refers to the multiplicative order of $k$ modulo $N$. Both are defined as the first positive integer $m$ such that $k^m\equiv 1\mod N$. Throughout this paper, $p$ will refer to primes. 

\subsection{Primitive Indexes in $k^u-h^u$}

\textbf{Proposition 1.} \hypertarget{P1}{}\emph{$N|k^u-h^u$ iff $P_N(k,h)|u$.}
\vspace{5mm}\\
Since
\begin{center}
$k^{P_N(k,h)}-h^{P_N(k,h)}\equiv 0\mod N\Longleftrightarrow k^{P_N(k,h)}\equiv h^{P_N(k,h)}\mod N$,
\end{center}
raising both $k$ and $h$ to a multiple of $P_N(k,h)$ will maintain the divisibility. In the other direction, since $N\nmid k^u-h^u$, where $u<P_N(k,h)$, then $k^u\not\equiv h^u\mod N$, and, for any positive integer $m$, 
\begin{center}
$k^{P_N(k,h)m+u}=k^{P_N(k,h)m}k^u\not\equiv h^{P_N(k,h)m}h^u=h^{P_N(k,h)m+u}\mod N$.
\end{center}
This proposition also shows that $N_1|N_2$ iff $P_{N_1}(k,h)|P_{N_2}(k,h)$ for all $k$ and $h$, providing the primitive indexes exist.
\vspace{5mm}\\
\textbf{Proposition 2.} \hypertarget{P2}{}\emph{If $gcd(N,k)=gcd(N,h)=1$, $P_N(k,h)|lcm(O_N(k),O_N(h))$ and $P_N(k,h)|\lambda(N)$.}\vspace{5mm}\\
Since $k^{\lambda(N)}\equiv h^{\lambda(N)}\equiv 1\mod N$ and $k^{lcm(O_N(k),O_N(h))}\equiv h^{lcm(O_N(k),O_N(h))}\equiv 1\mod N$, the result follows from the previous proposition. Additionally, $P_N(k,h)\neq O_N(k)$ and $P_N(k,h)\neq O_N(h)$ if $O_N(k)\neq O_N(h)$. It would be impossible for $k^{P_N(k,h)}\equiv h^{P_N(k,h)}\mod N$ if this were true. 
\vspace{5mm}\\
We will now present a modified version of a lifting the exponent [\hyperlink{3}{3}] proof to fit primitive indexes. 

\begin{theorem}\hypertarget{T1}
For $gcd(p,h)=1$, and for $k-h\neq 2(odd)$ and/or $p\neq 2$,
\begin{equation}
P_{p^a}(k,h)=pP_{p^{a-1}}(k,h),
\end{equation}\hypertarget{E1}{}
where $a-1\geq\nu_p(k^u-h^u)$, where $\nu_p(k^u-h^u)$ is the largest $i$ such that $u=P_{p^i}(k,h)=P_p(k,h)$.\\
For $gcd(p,h)=1$, and for $k-h=2(odd)$ and $p=2$, 
\begin{equation}
P_{2^a}(k,h)=2P_{2^{a-1}}(k,h),
\end{equation}\hypertarget{E2}{}
where $a-1\geq\nu_2(k^2-h^2)$, where $\nu_2(k^2-h^2)$ is the largest $i$ such that $2^i|k^2-h^2$. 
\end{theorem}
Take a prime $p$. If $p^a|k^u-h^u$, then for some positive integer $m$, $k^u-h^u=p^{a}m\Longrightarrow k^u=p^{a}m+h^u,$ where $p\nmid m$ and $a$ is the largest power of $p$ that divides $k^u-h^u$. Since $P_{p^a}(k,h)|P_{p^{a+1}}(k,h)$ by Proposition \hyperlink{P1}{1}, $P_{p^{a+1}}(k,h)$ is a multiple of $P_{p^a}(k,h)$. To search for this we do the following. For some integer $f\leq p$, we take the binomial expansion of $k^{uf}-h^{uf}=(p^am+h^u)^f-h^u$:
\begin{center}
$(\binom{f}{0}(p^am)^{f}+\binom{f}{1}(p^am)^{f-1}h^u+...+\binom{f}{f-2}(p^am)^{2}h^{u(f-2)}+\binom{f}{f-1}(p^am)^{1}h^{u(f-1)}+h^{uf})-h^{uf}$.
\end{center}
Then we cancel $h^{uf}$ and $-h^{uf}$ and factor out $p^a$: 
\begin{center}
$p^a(\binom{f}{0}p^{a(f-1)}m^{f}+\binom{f}{1}p^{a(f-2)}m^{f-1}h^u+...+\binom{f}{f-2}p^am^{2}h^{u(f-2)}+\binom{f}{f-1}mh^{u(f-1)}).$
\end{center}
If $\binom{f}{f-1}$ does not have a factor of $p$, then the expansion equals $pM+mh^{u(f-1)}\binom{f}{f-1}$ for some positive integer M. Recall that $p\nmid m$, and thus $p\nmid pM+mh^{u(f-1)}\binom{f}{f-1}$, as long as $gcd(p,h)\neq 1$. We thus specify $gcd(p,h)\neq 1$. This means that the only time when another factor of $p$ emerges is when the final term has a factor of $p$. For $f\leq p$, this only occurs when $f=p$: 
\begin{center}
$\binom{p}{p-1}mh^{u(f-1)}=\dfrac{p!}{(p-(p-1))!(p-1)!}mh^{u(f-1)}=pmh^{u(f-1)}.$
\end{center}
\emph{Case 1} ($a>1,$ $p$ odd): 
We can thus factor out exactly one more factor of $p$ for the first time at $k^{up}-h^{up}$, making $P_{p^{a+1}}(k,h)=up=pP_{p^a}(k)$, since the first time $p|k^u-h^u$ is $u=P_{p^a}(k,h)$.\\
\emph{Case 2} ($a=1$, $p$ odd): 
In only this case for odd numbers, the second last term is $\binom{p}{p-2}m^{2}h^{u(f-2)}$ after $p^2$ is factored out. But since $\binom{p}{p-2}m^{2}h^{u(f-2)}=\frac{p-1}{2}(p)(m^2)h^{u(f-2)}$, the second term still has a factor of $p$ and thus still yields $pM+mh^{u(f-1)}$ for the binomial as a whole, meaning there are no more factors of $p$.\\
\emph{Case 3} ($p=2$, $a=1$): 
After factoring out $2^2$ from $(2m+h^u)^2-h^{2u}$, we get: 
\begin{center}
$(2^2)(m^2+h^um)=(2^2)(m(m+h^u)).$
\end{center}
Since $2\nmid m$ or$ h^u$, $m$ and $h^u$ are odd and thus $m+h^u$ is even and has at least one more factor of $2$. This means that, if there is \emph{exactly} one factor of $2$ in $P_2(k,h)$, then $P_{2^x}(k,h)=2(P_2(k,h))$, where $x\geq 3$. This case is unique because $a=1$ as in Case $2$, and yet $2$ is even, meaning $2\nmid p-1=2-1$.\\
\emph{Case 4} ($p=2$, $a>1$):\\
After factoring out $2^{a+1}$, the final two terms of the binomial expansion are $\frac{2^{a+1}-2^a}{2}m^2h^{u(f-2)}+m=(2^{a}-2^{a-1})m^2h^{u(f-2)}+mh^{u(f-1)}=odd$, since neither $m$ nor $h$ can be even. It therefore cannot be divisible by $2$. This means that there are no more factors of $2$.\\
Equation (\hyperlink{E1}{1}) arises from the fact that either $k-h$ is odd and $2$ does not divide any term in the sequence, or $k-h$ is even such that $k-h\neq 2(odd)$, meaning $k-h=2(even)$, and thus the first term, $k-h$, is divided by at least $2^2$, and every power of $2$ greater than $1$ behaves normally. For equation (\hyperlink{E2}{2}), when $k-h=2(odd)$, it has exactly one factor of two, meaning that $k^2-h^2$ is divided by \emph{at least} $2^3$ by Case 3. Thus, the powers of $2$ will behave normally for sure only after the greatest power of $2$ that divides $k^2-h^2$. 
\vspace{5mm}\\
\textbf{Corollary.} \emph{For $gcd(p,h)=1$,}
\begin{center}
$P_{p^a}(k,h)=p^{a-\gamma_p(k,h)}P_p(k,h)$,
\end{center}
where 
\begin{center}
\[
  \gamma_p(k,h)=
  \begin{cases}
                                   \nu_p(k^u-h^u) & \text{if $k-h\neq 2(odd)$ and/or $p\neq 2$ and where $u=P_p(k,h)$;} \\
                                   \nu_2(k^2-h^2) & \text{if $k-h=2(odd)$ and $p=2$.} \\
  \end{cases}
\]
\end{center}
\textbf{Corollary.} \emph{
\begin{center}
$P_{nf^{a+1}}(k,h)=fP_n(k)$,
\end{center}
if $gcd(f,h)=1$ and the maximum power of $f$ that divides $k^{P_n(k,h)}-h^{P_n(k,h)}$ is $a$.}
\vspace{5mm}\\
\textbf{Corollary.} \emph{For $gcd(n,h_1)=gcd(n,h_2)=1$, iff $\gamma_n(k_1,h_1)=\gamma_n(k_2,h_2)$ and $P_n(k_1,h_1)=P_n(k_2,h_2)$, then $P_{n^a}(k_1,h_1)=P_{n^a}(k_2,h_2)$ for all positive integers $a$.}
\vspace{5mm}\\
\textbf{Proposition 3.} \hypertarget{P3}{}\emph{For $N=n_1n_2...$, where $gcd(N,h)=1$ and $gcd[n_i]_i=1$,}
\begin{center}
$P_N(k,h)=lcm[P_{n_i}(k,h)]_i$. 
\end{center}
From Proposition \hyperlink{P1}{1}, $n_i|k^u-h^u$ iff $P_{n_i}(k,h)|u$. All factors of coprime $n_i$ thus divide $k^u-h^u$ for the first time when $u=lcm[P_{n_i}(k,h)]$. 

\begin{theorem}\hypertarget{T2}{}
For $N=p_1^{p_1}p_2^{p_2}...$, $p_i$ prime, and $gcd(N,h)=1$,
\begin{center}
$P_N(k,h)=(\displaystyle\prod_i{n_i^{a_i-\gamma_{n_i}(k,h)})}lcm[P_{n_i}(k,h)]_i$.
\end{center}
\end{theorem}
We combine the previous properties to achieve this. If we let each $n_i=p_i$, $p_i$ prime, then $gcd[p_i]_i=1$, and we can find the order of any composite number coprime with the base $h$.
\vspace{5mm}\\
This reduces to 
\begin{center}
$O_N(k)=(\displaystyle\prod_i{n_i^{a_i-\gamma_{n_i}(k)})}lcm[O_{n_i}(k)]_i$,
\end{center}
for multiplicative order. 
\vspace{5mm}\\
We will now inspect certain sequences that yield the same primitive indexes for all numbers with specified prime factors. 
\begin{theorem}\hypertarget{T3}{} For any prime $p$, 
\begin{center}
$O_{p^a}(k+p^{\gamma_p(k)+1}m)=O_{p^a}(k),$
\end{center}
where $m\in\mathbb{N},m\neq 0$.\end{theorem}
For positive integers,
\begin{align*}
i^s-j^s&=(i-j)(i^{s-1}-i^{s-2}j+i^{s-3}j^2-...+i^2j^{s-3}-ij^{s-2}+j^{s-1}).
\end{align*}
This means, for any positive integer $m$, $k^n\equiv (k+pm)^n\mod p$ because $(k+pm)^n=k^n+((k+pm)^n-k^n)$, and $(k+pm)^n-k^n$ is divisible by $(k+pm)-k=pm$, as shown above. This means we are adding $0$ modulo $p$. It is therefore clear to see that $p$ will divide both $k^{O_p(k)}-1$ and $(k+pm)^{O_p(k)}-1$, and no exponents less than that. So $O_p(k)=O_p(k+pm)$. This holds true for all $p^a$ at $k+p^am$. Additionally, $k+p^am=k+p^b(p^{a-b}m')$.Therefore, at $p^a=p^{\gamma_p(k)+1}$, 
\begin{center}
$O_p(k+p^{\gamma_p(k)+1}m)=O_{p^2}(k+p^{\gamma_p(k)+1}m)=...=O_{p^{\gamma_p(k)}}(k+p^{\gamma_p(k)+1}m)=O_p(k),$
\end{center}
and
\begin{center}
$O_{p^{\gamma_p(k)+1}}(k+p^{\gamma_p(k)+1}m)=O_{p^{\gamma_p(k)+1}}(k)$,
\end{center} by the definition of $\gamma$. Since those two identities hold, we have proven that $O_{p^{\gamma_p(k)+1}}(k+p^{\gamma_p(k)+1}m)=pO_{p^{\gamma_p(k)+1}}(k)$ and all previous powers have the same order to that base, meaning  
\begin{center}
$\gamma_p(k+p^{\gamma_p(k)+1}m)=\gamma_p(k),$
\end{center} 
and $O_p(k+p^{\gamma_p(k)+1}m)=O_p(k)$. These two facts together are sufficient to prove the orders of $k$ modulo $p^a$ and $k+p^{\gamma_p(k)+1}m$ modulo $p^a$ are the same by the corollaries to Theorem \hyperlink{T1}{1}.
\vspace{5mm}\\
Let us denote the set of all positive integers made up solely of powers of $p_i$, for all chosen $i$, $i\in\mathbb{N}$, as follows: 
\begin{center}
$\eta[p_i]=\{p_1^{a_1}p_2^{a_2}...|a_i\in\mathbb{N}_0,p_i$ prime$\}$.
\end{center}
A number is said to be "in $\eta$" if it is one number made up of certain powers of some or all of the primes listed in the brackets. The period of any number, $N$, in $\eta$, base $k$ is defined as follows: 
\begin{center}
$\Lambda_N(k)=\Lambda_{\eta[p_i]}(k)=\displaystyle\prod_i{p_i^{\gamma_{p_i}(k)+1}}$.
\end{center}
\emph{If $p\nmid k^u-1$, then we define $\gamma_p(k)=0$ and if $k=1$, then $\Lambda_N(1)=0$.} Note that $\Lambda_N(k)$ is the same for every composite or prime number with the same unique prime factors, as they are all in $\eta$. The period of $\eta[p_i]$ base $k$ is defined as that of any and all $N$ in said $\eta$. \emph{Any} composite number or prime number $N$ in $\eta$ will behave the same in at least every base given as follows: 
\begin{theorem}\hypertarget{T4}{} For any number $N$ in $\eta[p_i]$,
\begin{center}
$O_{\eta[p_i]}(k+{\Lambda_{\eta[p_i]}(k)}m)=O_{\eta[p_i]}(k).$
\end{center} 
\end{theorem}
Since $O_{p^a}(k+p^{\gamma_p(k)+1}m)=O_{p^a}(k)$ from Theorem \hyperlink{T3}{3}, we know that $p_1^{\gamma_{p_1}(k)+1}m_1=p_2^{\gamma_{p_2}(k)+1}m_2=...$ will yield the above equality. This is the definition of $\Lambda$.

\begin{theorem}\hypertarget{T5}{} For $gcd(N,h)=1$ ,
\begin{center}$P_{\eta[p_i]}(k+\Lambda_{\eta[p_i]}(k)m,h+\Lambda_{\eta[p_i]}(h)m)=P_{\eta[p_i]}(k,h)$.
\end{center} 
\end{theorem}
$(j+\Lambda_{\eta[p_i]}(j)m)^u\equiv j^u\mod N$ by Theorem \hyperlink{T4}{4}. This means that the residues for all exponents are the same for bases $k$ and $k+\Lambda_{\eta[p_i]}(k)m$ and $h$ and $h+\Lambda_{\eta[p_i]}(h)m$. The result follows after taking into account that when $p\nmid k^u-1$, $\gamma_p(k)=0$ and if $k=1$, then $\Lambda_N(1)=0$. 
\vspace{5mm}\\
\textbf{Remark.} Note that, for $gcd(k,h)=1$, no number $gcd(n,h)\neq 1$ and $gcd(n,k)\neq 1$ will divide $k^u-h^u$, and thus every number that does divide $k^u-h^u$ will obey Theorem \hyperlink{T1}{1} and Proposition \hyperlink{P2}{2}.

\subsection{Primitive Indexes in $k^u+h^u$}
We operate under the condition that neither $k^u+h^u$ nor $k^u-h^u$ are equal to $0$ for all $u$.  
\begin{theorem}\hypertarget{T6}{}
For $N>2$, $N|k^u+h^u$ iff $P^N(k,h)|u$ and 
\begin{center}
$\dfrac{u}{P^N(k,h)}$,
\end{center}
is odd.
\end{theorem}
When $N|k^{P^N(k,h)}+h^{P^N(k,h)}$, we have 
\begin{center}
$k^{P^N(k,h)}+h^{P^N(k,h)}\equiv 0\mod N\Longleftrightarrow k^{P^N(k,h)}\equiv -h^{P^N(k,h)}\mod N$. 
\end{center}
Raising both $k$ and $h$ to the same odd multiple of $P^N(k,h)$ will preserve this congruence, even multiples will not. In the other direction, for all $y<P^N(k,h)$, we have $k^y\not\equiv h^y\mod N$. This is proven as follows: we see that for all $N$ that divide some $k^u+h^u$, 
\begin{center}
$k^{P^N(k,h)}\equiv -h^{P^N(k,h)}\mod N\Longrightarrow k^{2P^N(k,h)}\equiv h^{2P^N(k,h)}\mod N,$
\end{center}
and by Proposition \hyperlink{P1}{1}, $P_N(k,h)|2P^N(k,h)$. But since $P_N(k,h)\nmid P^N(k,h)$ for all $N>2$ by Proposition \hyperlink{P1}{1}, for all such $N$ we know $P_N(k,h)=2P^N(k,h)$ or $P_N(k,h)=2$. But if $P_N(k,h)=2$, then because $k^{2}-h^{2}=(k-h)(k+h)$ and $N\nmid k-h$, $N|k+h$ and by definition $P^N(k,h)=1$. In this case, it is still true that $2=2P^N(k,h)$. Thus, for all $y<P^N(k,h)<P_N(k,h)$, we have
\begin{center}
$k^{P^N(k,h)m+y}=k^{P^N(k,h)m}k^y\not\equiv -h^{P^N(k,h)m}h^y=-h^{P^N(k,h)m+y}\mod N$,
\end{center}
for any positive integer $m$. For $N=2$, $N$ either divides every term or no term. 

\begin{theorem}\hypertarget{T7}{}
$2P^N(k,h)=P_N(k,h)$, where $N=p^a$ and $p$ has an even primitive index in $k^u-h^u$ and/or where $N>2$ divides some term in $k^u+h^u$.
\end{theorem}
We proved that for all $N$ which divide some $k^u+h^u$, the equality holds in the theorem above in Theorem \hyperlink{T6}{6}. In the other direction, for an \emph{odd} prime power $p^a$ with an even primitive index in $k^u+h^u$,
\begin{align*}
k^{P_{p^a}(k,h)}-h^{P_{p^a}(k,h)}&= (k^{\frac{P_{p^a}(k,h)}{2}}-h^{\frac{P_{p^a}(k,h)}{2}})(k^{\frac{P_{p^a}(k,h)}{2}}+h^{\frac{P_{p^a}(k,h)}{2}})\\
&\equiv 0\mod p^a,
\end{align*}
implies that $p^a|(k^{\frac{P_{p^a}(k,h)}{2}}+h^{\frac{P_{p^a}(k,h)}{2}})$ since $p\nmid (k^{\frac{P_{p^a}(k,h)}{2}}-h^{\frac{P_{p^a}(k,h)}{2}})$ by Theorem \hyperlink{T1}{1} and Proposition \hyperlink{P1}{1}. This is because any odd $p^a$ with an even primitive index we will have $y=\frac{P_{p^a}(k,h)}{2}=p^{a-\gamma_{p^a}(k,h)}\frac{P_p(k,h)}{2}$ and thus $P_p(k,h)\nmid y$. Additionally, no smaller exponent $u$ will allow $k^u\equiv -h^u\mod {p^a}$ because then $p^a|k^{2u}-h^{2u}$ for $2u<P_{p^a}(k)$ by the forward direction of this theorem proven in Theorem \hyperlink{T6}{6}, contradicting Proposition \hyperlink{P1}{1}. For all odd primes, $P_{p}(k,h)$ and $P_{p^a}(k,h)$ are both either odd or even. Therefore, $P_{p}(k,h)$ being even implies $P^{p^a}(k,h)=\frac{1}{2}P_{p^a}(k,h)$. 

\begin{theorem}\hypertarget{T8}{} For odd $N=p_1^{a_1}p_2^{a_2}...$, iff each $P^{p_i^{a_i}}(k,h)$ is odd or the smallest primitive index in the $lcm$ is even and every other primitive index in the $lcm$ is an odd multiple of it, 
\begin{center}
$P^N(k,h)=lcm[P^{p_i^{a_i}}(k,h)]$. 
\end{center}
\end{theorem}
By Theorem \hyperlink{T6}{6} and the logic from Proposition \hyperlink{P3}{3}, the result follows when $P^{p_i^{a_i}}(k,h)$ is odd and/or for odd multiples of the smallest primitive index. \emph{For any other primitive indexes, there can never be a multiple of them such that dividing by any one of them always yields an odd, meaning $N\nmid k^u+h^u$ for any $u$ if the conditions specified in the theorem are not met.} 
\vspace{5mm}\\
For odd $p$, $P^{p^a}(k,h)$ may be calculated as in Theorem \hyperlink{T1}{1}, as shown by Theorem \hyperlink{T6}{6}. 
\vspace{5mm}\\
We can finalize the discussion of lifting the exponent in [\hyperlink{3}{3}] by deriving a formula for $\nu_2(k^n+h^n)$, as I have not found it  discussed elsewhere. Recall that $\nu_2(k^n+h^n)$ refers to the highest power of $2$ that divides $k^n+h^n$. Also note that $P_2(k,h)=P^2(k,h)=1$ if it exists, and $2$ either divides every term in $k^u+h^u$ and $k^u-h^u$ or none. 
\begin{theorem}\hypertarget{T9}{} 
Iff $2|k+h$,
\[
  \nu_2(k^n+h^n) =
  \begin{cases}
                                   \nu_2(k+h) & \text{if $n$ is odd;} \\
                                   1 & \text{if $n$ is even.} \\
  \end{cases}
\]
\end{theorem}
If $2\nmid k+h$, then clearly $2$ divides no $k^n+h^n$. For $\nu_2(k+h)=a$, we know $2^{a-x}\nmid k^n+h^n$ for even $n$ and $a-x>1$ by Theorem \hyperlink{T6}{6}. Since if $2|k+h$, $2|k^n+h^n$, it must be true that $\nu_2(k^n+h^n)=1$ for even $n$. For odd $n$, we have that $2^a|k^n+h^n$ by Theorem \hyperlink{T6}{6}. If $2^{a+x}|k^n+h^n$, then the first index it does this at yields $P_{2^{a+x}}(k,h)=2P^{2^{a+x}}(k,h)=2(odd), odd\neq 1$ again by Theorem \hyperlink{T7}{7}. But since for any $y$, $P_{2^y}(k,h)$ is a power of $2$ by Theorem \hyperlink{T1}{1}, this is a contradiction. 
\vspace{5mm}\\
By combining Theorem \hyperlink{T8}{8} and Theorem \hyperlink{T9}{9}, one can easily see whether any even number $N$ has a primitive index in $k^u+h^u$ and furthermore calculate the primitive index if it does. 
\section{Connection to Primitive Roots}

We know that the multiplicative group of integers modulo $n=2,4,p^a,$ or $2p^a$, $(\mathbb{Z}/n\mathbb{Z})^{\times}$, has $|(\mathbb{Z}/n\mathbb{Z})^{\times}|=\phi(n)$. Its generators can be expressed as ${k}$, where $O_n(k)=\phi(n)$. This means that 
\begin{center}
$\{k^0,k^1,...,k^{\phi(n)-1}\}$,
\end{center}
represents a reduced residue system modulo $n$. At least one $k$ yields this iff $n=2,4,p^a,$ or $2p^a$. This information is from [\hyperlink{1}{1}]. 
\vspace{5mm}\\
\textbf{Proposition 4.}\hypertarget{P4}{}
\emph{For $0\leq u<\phi(n)$, $k^u\equiv h^u\mod n$ exactly $t$ times, where $P_n(k,h)t=\phi(n)$.}
\vspace{5mm}\\
By Proposition \hyperlink{P1}{1}, for $u<\phi(n)$, $k^u\equiv h^u\mod n$ exactly at $u=0,P_n(k,h),2P_n(k,h),...,(t-1)P_n(k,h)$. The next equivalence is not less than $\phi(n)$. 
\vspace{5mm}\\
This proposition says that $P_n(k,h)=\phi(n)$ means both 
\begin{center}
$\{k^0,k^1,...,k^{\phi(n)-1}\}$,
\end{center}
and
\begin{center}
$\{h^0,h^1,...,h^{\phi(n)-1}\}$,
\end{center}
are reduced residue systems modulo $N$ and furthermore that $k^0\equiv h^0\mod n$ and every other equivalent power of $k$ and $h$ in the reduced residue system are distinct modulo $n$. 
\vspace{5mm}\\
\textbf{Proposition 5.} \emph{
If $k(h)^{\phi(n)-P_n(k,h)}$ is a generator of the reduced residue system of $n$, then so are $k$, $h$, and $h(k)^{\phi(n)-P_n(k,h)}$.}
\vspace{5mm}\\
First note that, by Proposition \hyperlink{P4}{4}, $O_n(k)=O_n(h)=\phi(n)$, and both $h$ and $k$ are also generators of the reduced residue system. This means the inverse of $a^{b}\mod n$ is $a^{\phi(n)-b}\mod n$. Therefore, 
\begin{center}
$k^{P_n(k,h)}-h^{P_n(k,h)}\equiv 0\mod n\Longleftrightarrow h^{P_n(k,h)}(k^{P_n(k,h)})^{\phi(n)-P_n(k,h)}\equiv k^{P_n(k,h)}(h^{P_n(k,h)})^{\phi(n)-P_n(k,h)}\equiv 1\mod n$. 
\end{center}
Additionally, no exponent $u<P_n(k,h)$ will yield either equivalence. This means that $P_n(k,h)=O_n(k(h)^{\phi(n)-P_n(k,h)})=O_n(h(k)^{\phi(n)-P_n(k,h)})$. Therefore, if $P_n(k,h)=\phi(n)$, then so do those two orders. 

\section{Cyclotomic Polynomials, Zsigmondy Numbers, and Primoverization}
A homogenized cyclotomic polynomial is defined as follows:
\begin{center}
$\Phi_n(k,h)=\dfrac{k^n-h^n}{\displaystyle\prod_{d|n} {\Phi_{d}(k,h)}}$.
\end{center}
Straightforwardly, the highest power of all nonprimitive factors of $k^n-h^n$ that divide some $k^d-h^d$ for $d|n$ are divided out via this definition. 
\vspace{5mm}\\
\textbf{Definition: Zsigmondy Numbers.}
A Zsigmondy number,
\begin{center}
$\mathcal{Z}(n,k,h)$,
\end{center}
is the product of the primitive prime divisors of the term $k^n-h^n$. We shall denote 
\begin{center}
$\zeta(n,k,h)$,
\end{center}
to be the product of the primitive prime divisors of the term $k^n+h^n$. These are also termed Zsigmondy numbers. 
\vspace{5mm}\\
One can see how Zsigmondy numbers and primitive indexes are closely related. Without so much as invoking the properties of cyclotomic polynomials, we can find Zsigmondy numbers explicitly. 

\begin{theorem}\hypertarget{T10}{}
For $n\neq 2$ and $gcd(k,h)=1$ and $z\in\mathbb{N},z\neq 0$, 
\[
  \mathcal{Z}(n,k,h) =
  \begin{cases}
                                   \Phi_n(k,h) & \text{if $n\neq p^z(P_p(k,h))$;} \\
                                   \dfrac{\Phi_n(k,h)}{p} & \text{if $n= p^z(P_p(k,h))$.} \\
  \end{cases}
\]
\end{theorem}
As noted above, the highest power of all nonprimitive factors of $k^n-h^n$ which divide some $k^d-h^d$ for $d|n$ are divided out via the definition of $\Phi_n(k,h)$. This means the only nonprimitive prime factors of $\Phi_n(k,h)$ are new multiplicities of previous prime factors introduced in $k^n-h^n$. From Theorem \hyperlink{T2}{2}, the new multiplicities of $N$ can only be present at $n=N^zq$, where $N|k^q-h^q$. Additionally from the theorem, the highest power of $N|\Phi_n(k,h)$, where $N$ is not primitive and $n>2$, is $1$. 
\vspace{5mm}\\
\textbf{Lemma.} \emph{If $q\neq P_N(k,h)$, then $N\nmid\Phi_n(k,h)$.}
\vspace{5mm}\\
Assume $q\neq P_N(k,h)$. By Proposition \hyperlink{P1}{1}, $P_N(k,h)|q$. By Theorem \hyperlink{T2}{2}, we can see 
$P_{N^a}(k,h)=N^{a-\gamma_N(k,h)}P_N(k,h)|N^zq$. This would be divided out. 
\vspace{5mm}\\
\textbf{Lemma.} \emph{If $N$ is composite, then $N\nmid\Phi_n(k,h)$.}
\vspace{5mm}\\
Assume $N$ is composite. Since this must be a new multiplicity of $N$, we know $P_N(k,h)<n$. By Property \hyperlink{P1}{1}, we see that $P_{N^a(f_i)}(k,h)|P_{N^{a+1}}(k,h)=n$ for every proper factor $f_i$ of $N$ and for $a+1$ the largest power of the chosen $f_i$ that divides $k^n-h^n$. By Theorem \hyperlink{T2}{2} and because $n>2$, we have that $P_{N^a(f_i)}(k,h)<P_{N^{a+1}}(k,h)$, and thus $f_i\nmid \Phi_n(k,h)$. All factors $f_i$ of $N$ would be divided out.
\vspace{5mm}\\
It is clear that $n=p^zP_p(k,h)$ if it possible that $p|\Phi_n(k,h)$ and $p$ is not primitive ($p$ prime). By Theorem \hyperlink{T2}{2}, the maximum power of $p$ that divides $\Phi_n(k,h)$ is $1$. In the other direction, if $n= p^z(P_p(k,h))$, then $n=P_{p^a}(k,h)$, where $a=\nu_p(k^n-h^n)$, by Theorem \hyperlink{T2}{2}, and thus there must be exactly one factor of $p$ that divides $\Phi_{p^zP_p(k,h)}(k,h)$ and is not primitive. 
\vspace{5mm}\\
This reduces to 
\[
  \mathcal{Z}(n,k,1) =
  \begin{cases}
                                   \Phi_n(k) & \text{if $n\neq p^z(O_p(k))$;} \\
                                   \dfrac{\Phi_n(k)}{p} & \text{if $n= p^z(O_p(k))$.} \\
  \end{cases}
\]
This definition is preferred because it allows one to classify certain classes of numbers as Zsigmondy numbers without even calculating $\Phi_n(k,h)$. 
\vspace{5mm}\\
\textbf{Proposition 6.}
\emph{For $n=2$, 
\begin{center}
$\mathcal{Z}(2,k,h)=\dfrac{\Phi_2(k,h)}{2^{\nu_2(k^2-h^2)}}$.
\end{center}}
This follows directly from Theorem \hyperlink{T1}{1} and the reasoning expressed in Theorem \hyperlink{T10}{10}. 
\vspace{5mm}\\
We have not only proved the forward condition that $n=p^zq$, but additionally that $q=P_p(k,h)$ and the biconditionality of the statement. We will now explore how this affects Shevelev's exploration of overpseudoprimes. From Shevelev [\hyperlink{2}{2}, (2.4)], we see that all primes satisfy 
\begin{center}
$p=r_k(p)O_p(k)$,
\end{center}
for any base $k$, where $r_k(p)$ is the number of distinct cyclotomic cosets of $k$ modulo $p$. A composite number is termed \emph{overpseudoprime} base $k$ if $n=r_k(n)O_n(k)$.
\vspace{5mm}\\
\textbf{Definition: Primover.}
A number $n$ is said to be primover base $k$ if it is either a composite overpseudoprime base $k$ or prime. 
\vspace{5mm}\\
From Shevelev [\hyperlink{2}{2}, Th. 12], a necessary and sufficient condition for $n$ to be primover is as follows: 
\vspace{5mm}\\
\textbf{Lemma.} \emph{If $gcd(n,b)=1$, then $n$ is primover iff $O_n(k)=O_d(k)$ for each divisor $d>1$ of $n$.}
\vspace{3mm}\\
In his paper [\hyperlink{2}{2}, Th. 16, Th. 17], he details how every composite overpseudoprime base $k$ is a strong pseudoprime and superpseudoprime to the same base. We strengthen and expand upon some of the theorems found in [\hyperlink{2}{2}], making them biconditional and removing the condition "if $p$ does not divide $\Phi_n(k)$" (if $gcd(p,\Phi_n(k))=1$). We also introduce new examples of primover numbers evident from Theorem \hyperlink{T10}{10}. Note that all Zsigmondy numbers in the sequence $k^u-1$ are primover due to the above lemma. 
\vspace{5mm}\\
\textbf{Proposition 7.} \emph{Generalized repunits,
\begin{center}
$\mathcal{Z}(n,k,1)=\dfrac{k^n-1}{k-1}$,
\end{center}
are primover base $2$ iff $n$ is prime.}
\vspace{5mm}\\
\textbf{Proposition 8.} \emph{Numbers of the form\begin{center}
$\mathcal{Z}(pq,k,1)=\dfrac{(k-1)(k^{pq}-1)}{(k^p-1)(k^q-1)},$
\end{center}
are primover base $k$. This identity holds iff both $p$ and $q$ are distinct primes, where $p\neq O_q(k)$ and $q\neq O_p(k)$.}
\vspace{5mm}\\
\textbf{Proposition 9.} \emph{Numbers of the form\begin{center}
$\mathcal{Z}(pO_p(k),k,1)=\dfrac{(k-1)(k^{pO_p(k)}-1)}{p(k^p-1)(k^{O_p(k)}-1)}$
\end{center}
are primover base $k$ iff $p$ is prime and $O_p(k)$ is prime.}
\vspace{5mm}\\
\textbf{Proposition 10.} \emph{Numbers of the form 
\begin{center}
$\mathcal{Z}(n^a,k,1)=\dfrac{k^{n^a}-1}{k^{n^{a-1}}-1}$,
\end{center}
are primover base $k$ iff $n$ is prime.} 
\vspace{5mm}\\
We further the study of Zsigmondy by inspecting $k^u+h^u$. We introduce an important theorem regarding the relationship between $k^u+h^u$ and $k^u-h^u$.
\begin{theorem}
For $n>1$ and/or $k+h$ odd, 
\begin{center}
$\zeta({n},k,h)=\mathcal{Z}({2n},k,h)$.
\end{center}
\end{theorem}
From Theorem \hyperlink{T7}{7}, if an odd $p$ divides some $k^u+h^u$, then $P_{p^a}(k,h)=2P^{p^a}(k,h)$. If an odd $p$ has an even primitive index in $k^u-h^u$, then by Theorem \hyperlink{T7}{7}, $P_{p^a}(k,h)=2P^{p^a}(k,h)$. By the definition of primitive prime divisors, the result follows when $n>1$, since $P_2(k,h)=P^2(k,h)=1$. If $k+h$ is odd then so will every term be, and $\zeta(1,k,h)=\mathcal{Z}(1,k,h)$. If not, then since $P_2(k,h)=P^2(k,h)=1$, $2^{a-b}\zeta(1,k,h)=\mathcal{Z}(1,k,h)$, where $a$ is the highest power of $2$ that divides $k-h$ and $b$ is the highest power of $2$ that divides $k+h$. Thus, all $\zeta(n,k,1)$ are primover base $k$ for $n>1$ and/or $k+h$ odd. For $n=1$ and $k+1$ even, $\zeta(1,k,1)$ is primover iff $a>b$.
\vspace{5mm}\\
\textbf{Corollary.} \emph{The only odd $N$ that divide some $k^u+h^u$ are those that are made up entirely of primes with even primitive indexes in $k^u-h^u$ and that satisfy the conditions specified in Theorem \hyperlink{T8}{8}.} 
\vspace{5mm}\\
\textbf{Corollary.} \emph{For $k+u$ even, the only even $2^aN$, $N$ odd, that divide some $k^u+h^u$ are those where $N$ obeys the above corollary and $a=\nu_2(k+h)$ if $P^N(k,h)$ is odd and $a=1$ if $P^N(k,h)$ is even.}
\vspace{5mm}\\
This is proven via Theorem \hyperlink{9}{9}. 
\vspace{5mm}\\
\textbf{Proposition 11.} \emph{
All Wagstaff numbers with prime exponents, 
\begin{center}
$\zeta(p,2,1)=\dfrac{2^p+1}{3}$,
\end{center}
are primover base $2$.}

\newcommand{\noopsort}[1]{} \newcommand{\printfirst}[2]{#1}
  \newcommand{\singleletter}[1]{#1} \newcommand{\switchargs}[2]{#2#1}

\end{document}